\documentclass[11pt]{article}

\usepackage{amsmath,amsfonts,amscd, amssymb,theorem}

\newcommand{\Span}[1]{\left<#1\right>}
\newcommand{\doubleheaddownarrow}{\big\downarrow\kern-3.325mm\downarrow}
\newcommand\Oh{\mathcal O}
\newcommand\Z{\mathbb Z}
\newcommand\Q{\mathbb Q}
\newcommand\C{\mathbb C}

\newcommand\om{\omega}

\newcommand{\iso}{\cong}

\newcommand{\Ann}{\operatorname{Ann}}
\newcommand{\commandeta}{\eta}

\newcommand\Proj{\operatorname{Proj}}

\newcommand\Hom{\operatorname{Hom}}

\newcommand{\skewentry}{ \kern -1cm -\mathrm{sym} \kern -1cm}
\newcommand\Spec{\operatorname{Spec}}

 \renewcommand{\P}{\operatorname {\mathbb P}}

 \newcommand\A{\operatorname {\mathbb A}}

 \newtheorem{theorem}{Theorem}[section]
 \newtheorem{lemma}[theorem]{Lemma}
 \newtheorem{prop}[theorem]{Proposition}
 \newtheorem{cor}[theorem]{Corollary}
 {
 \theorembodyfont{\rmfamily}
 \newtheorem{defn}[theorem]{Definition}

 \newtheorem{rem}[theorem]{Remark}
 
 \newtheorem{notation}[theorem]{Notation}

 }
 \newenvironment{pf}{\paragraph{Proof}}{\par\medskip}
 
 \newcommand{\qed}{\ifhmode\unskip\nobreak\fi\quad\ensuremath\square}
\newcommand{\QED}{\ifhmode\unskip\nobreak\fi\quad\ensuremath{\mathrm{QED}}}

\numberwithin{equation}{section}

\renewcommand{\labelenumi}{(\arabic{enumi})}

\newcommand{\set}[1]{\left\{ #1 \right \} }
\newcommand{\dd}[2]{\frac{\partial{#1}}{\partial{#2}}}

\title{A construction of numerical Campedelli surfaces with torsion  $\Z/6$ }

\author{Jorge Neves and Stavros Argyrios Papadakis}
\date {}

\begin{document}
\maketitle

\begin {abstract}
We produce a family of numerical Campedelli surfaces with  $\Z/6$ torsion  by
constructing the canonical ring of the \'{e}tale six to one cover
using serial unprojection. In Section~\ref{sec!algebrapart} we develop
the necessary algebraic machinery. Section~\ref{sec!campedelli} contains
the  numerical Campedelli surface construction, while
Section~\ref{sec!finalcom} contains remarks and open questions.
\end{abstract}

\section  {Introduction}  
A numerical Campedelli surface is a smooth minimal surface of general type over
the complex numbers with $K^2 =2$ and $q=p_g=0$. It is known that the algebraic
fundamental group $\pi_1^{\mathrm{alg}}$ of such surface is finite, of order at
most 9 (cf.~\cite{BPHV} Chap.~VII.10). Two recent papers about numerical Campedelli
surfaces  are \cite{MP} and \cite{LP}, the first classifies the case where the
algebraic fundamental group has order exactly $9$, while the second gives
simply connected examples.

In the present work we give a construction of
numerical Campedelli surfaces with algebraic fundamental group  equal to $\Z/6$.  To our
knowledge, there were no such examples previously known, and it settles the
existence question for numerical Campedelli surfaces with algebraic fundamental group of
order $6$, since by \cite {Na} there are no numerical  Campedelli surfaces with
algebraic fundamental group  equal to the symmetric group of order six.

Our approach is to construct, using serial unprojection of type Kustin--Miller,
the canonical ring of  the \'{e}tale  six to one cover together with a suitable
basepoint free action of  $\Z/6$. The cover is a regular canonical surface with $p_g=5$,
and $K^2=12$, canonically embedded in $\P(1^5,2^3)$.

In Section~\ref{sec!algebrapart} we define, for $n \geq 2$, what we call the
generic $\binom {n} {2}$ Pfaffians ideal
(Definition~\ref{dfn!binomfaffiansideal}) and prove that it is Gorenstein of
codimension equal to $n+1$ (Theorem~\ref{thm!mainalgebrathm}). A
special case of the construction is due to Frantzen (\cite {Fr}  Section~2.4).

 In Section~\ref{sec!campedelli} we apply the results of the previous sections in
the case of $n=4$ to our specific geometric situation. The main results are
Theorems~\ref{theorem!smoothnessandirreducibility} and 
\ref{theorem!existenceCampedelli} where we settle the existence of a 
nonsingular regular surface with $p_g=5$ and $K^2=12$ endowed with a $\Z/6$
basepoint free action. Finally, Section~\ref{sec!finalcom} contains some remarks and open
questions.

The way in which we have arrived to the family constructed in this
article is strongly influenced by the general theory in \cite{R1}.
More precisely, one assumes that an hypothetical \'{e}tale
six to one cover of a numerical Campedelli surface is a quadratic
section of an anticanonically embedded Fano threefold $V$, as in many
other examples of surfaces of general type. Then, standard numerical
Hilbert series calculations (cf.~\cite{R1} Section~3) lead
to the expectation that the anticanonical model of $V$ is a
codimension $5$ projectively Gorenstein subscheme of $\P(1^5,2^4)$
with a certain Hilbert series. Combining this with the knowledge
of how the Hilbert series changes during unprojection (or,
alternatively, and more easily, read-off the result directly from
\cite{Br}) one realizes that $V$ can be realized as the result
of a series of four unprojections of Kustin--Miller type, starting
from a degree $4$ hypersurface in $\P(1^5)$. Hence, starting from a
degree $4$ hypersurface in $\P(1^5)$ and unprojecting an arrangement of
$4$ codimension $1$ loci one could obtain a $3$-fold $V$ in $\P(1^5,2^4)$
with the right Hilbert series. Then, taking a suitable member of
$|-2K_V|$ we would obtain the \'{e}tale six to one cover of a
numerical Campedelli.

However, to set up a free $\Z/6$ action, motivated by empirical
evidence showing that unprojection is best calculated in a general
framework, we were driven to the unprojection of a general set of $4$
linear subspaces of dimension $5$ in a degee $4$ $6$-fold hypersurface 
in $\P(1^8)$.
Our main motivation was that we could then assume that these loci were
defined by $x_1=x_2=0$, $x_3=x_4=0$, etc. After the unprojection of
these subspaces we obtained a $6$-fold in $\P(1^8,2^4)$ and then taking
$3$ linear sections and $1$ quadratic section we constructed a family
of surfaces with $p_g =5, q=0, K^2=12$ (cf. Remark~\ref{rem!forgetaction} below).
Afterwards, calculations with characters and $G$-Hilbert series helped
us to discover a suitable subfamily endowed with a good $\Z/6$ action.

We believe that a similar approach could be useful to other
situations as well, compare Remark~\ref{remark!no335} below.

The expectation that under mild conditions unprojections commute 
(cf.~Remark~\ref{rem!commutativityofunprojections} below) and bearing in mind 
the previously done calculations of
\cite{P2} and \cite{Fr}~Section~2.4 we got to the
the generic $\binom {n} {2}$ Pfaffians  ideal format for $n=4$. 
We then discovered that the arguments for the Gorensteiness of the format 
for general $n$ were a rather straightforward generalisation of those
needed in the special case.

\medskip

\paragraph{Acknowledgments} The authors wish to thank Miles Reid and Margarida Mendes
Lopes for important discussions and suggestions which have impro\-ved the
presentation of the paper. The authors are grateful for the financial support
of Funda\c{c}\~ao Gulbenkian. The second author was supported by the Portuguese 
Funda\c{c}\~ao para a Ci\^encia e a Tecno\-lo\-gia through Grant SFRH/BPD/22846/2005 
of POCI2010/FEDER.

\section {The generic $\binom {n} {2}$ Pfaffians ideal }  \label {sec!algebrapart}

\begin{notation}   \label {not!bignotation}
Let us make the following notation.

\begin{enumerate}
\item Let $n \geq 2$. Let $A_0 = \Z [x_1, \dots ,x_n, z_1, \dots ,z_n, r_{d_1\cdots d_n} ]$,
be the polynomial ring over the integers in $n$ variables $x_1, \dots ,x_n$, $n$ variables 
$z_1, \dots ,z_n$ and $2^n$ variables 
$r_{d_1 \cdots d_n}$, indexed by $(d_1, \dots, d_n) \in \{0,1 \}^{n}$.

\item Define the polynomial algebra extensions
$A_0 \subset  A_1 \subset A_2 \subset \dots \subset  A_n$ by setting
inductively $A_{i} = A_{i-1} [y_i]$ for $i =1, \dots ,n$.

\item Make these rings graded by setting the degree of $x_i$, $z_i$ and
$r_{d_1 \cdots d_n}$ equal to $1$, for all $i=1,\dots, n$ and
$(d_1, \dots ,d_n)\in\set{0,1}^n$ and by setting the degree of $y_i$ equal to
$n-1$, for all $i=1,\dots,n$.

\item Consider the degree $n+1$ homogeneous polynomial defined by
\[
   Q = \sum r_{d_1 \cdots d_n} a_{1,d_1} \cdots a_{n,d_n} \in A_0,
\]
where the summation is for  $(d_1, \dots ,d_n) \in \{ 0,1 \}^{n}$, and, by
definition,  $a_{i,d_i}$ is equal to $x_i$ if $d_i =0 $ and equal to $z_i$ if
$d_i =1$.

\item  \label{item!105}
For each $1\leq i< j\leq n$, let
\[
\begin{array}{c}
Q_{ij}^{xx}=\dd{^2Q}{x_i\partial x_j},\quad Q_{ij}^{xz}=
\dd{^2Q}{x_i\partial z_j},\quad  Q_{ij}^{zx}= \dd{^2Q}{z_i\partial x_j}\quad
\text{and}\quad Q_{ij}^{zz}= \dd{^2Q}{z_i\partial z_j}.
\end{array}
\] Then each of the
$Q_{ij}^{ab}$, where $a,b\in\set{x,y}$, is homogeneous of degree $n-1$ and
clearly,
\[
   Q = x_ix_jQ_{ij}^{xx}+ x_iz_jQ_{ij}^{xz}+z_ix_jQ_{ij}^{zx}+ z_iz_j Q_{ij}^{zz}.
\]

\item \label{item!101}
 For each $1\leq i< j\leq n$ let
\[
    M_{ij} =  \begin{pmatrix}
   0   & x_i & z_i & -x_j & -z_j \\
       &   0 &  y_j & Q^{zz}_{ij} & -Q^{zx}_{ij}\\
       &     &  0   & -Q^{xz}_{ij} & Q^{xx}_{ij} \\
       &  \text{-sym}   &       & 0   & -y_i \\
       &    &      &   &   0
    \end{pmatrix}
\]
be a skew--symmetric $5 \times 5$ matrix with entries in $A_n$. The $5$ submaximal Pfaffians\footnote 
{For a discussion about the Pfaffians  of a skew--symmetric matrix  see, 
for example, \cite {BH} Section~3.4.}
of this matrix are:
\begin{equation}\label{eqn!pfaffians}
\renewcommand{\arraystretch}{1.4}
\begin{array}{c}
y_iy_j- Q_{ij}^{xz}Q_{ij}^{zx}+Q_{ij}^{xx}Q_{ij}^{zz}, \\
x_iy_i + (x_jQ_{ij}^{zx}+z_jQ_{ij}^{zz}),  \quad   \quad  z_iy_i -
(x_jQ_{ij}^{xx}+z_jQ_{ij}^{xz}), \\
 x_jy_j + (x_iQ_{ij}^{xz}+z_iQ_{ij}^{zz}), \quad  \quad z_jy_j -
 (x_iQ_{ij}^{xx}+z_iQ_{ij}^{zx});
\end{array}
\end{equation}
which are all homogeneous elements of $A_n$. Fixing $1\leq i\leq n$ and varying $j$
we see that several Pfaffians involve the monomial $x_iy_i$. Notice, however,
that $x_jQ_{ij}^{zx}+z_jQ_{ij}^{zz}=\dd{Q}{z_i}$, which does not depend on $j$.
Hence, in the set of Pfaffians of all possible $M_{ij}$, there is only one
polynomial in which the monomial $x_iy_i$ occurs. A similar observation
applies to the Pfaffians in which the monomial $z_iy_i$ occurs.

\item For each $0 \leq p \leq n$ we define an homogeneous ideal
$I_p\subset A_p$ by:
\begin{enumerate}
\item $I_0 =(Q) \subset A_0$;
\item $I_1= \left (x_1y_1+\dd{Q}{z_1},z_1y_1-\dd{Q}{x_1} \right )\subset A_1$,
the ideal of $A_1$ generated by the two Pfaffians of $M_{12}$ which involve
$x_1y_1$ and $z_1y_1$;
\item For $2 \leq p \leq n$, $I_p\subset A_p$ is the ideal of $A_p$
generated by the union of all the submaximal Pfaffians of all
matrices $M_{ij}$ for $1 \leq i < j\leq p$.
\end{enumerate}

\item \label{item!102}
In the set of Pfaffians of all possible $M_{ab}$, $1 \leq a < b \leq n$, denote by
\begin{enumerate}
\item    $e^{xy}_i$ the polynomial in which $x_iy_i$ occurs with coefficient $1$;
\item  $e^{zy}_i$ the polynomial in which $z_iy_i$ occurs with coefficient $1$;
\item  $e^y_{ij}$ the polynomial in which $y_iy_j$ occurs with coefficient $1$.
\end{enumerate}
for all $1\leq i \leq n$ in (a) and (b), and  
    $1\leq i < j \leq n$  in  (c). In particular,
\[
I_p = \left ( \set{e^y_{ij}, e^{xy}_t, e^{zy}_t : 1 \leq i < j \leq p \text{
and } 1 \leq t \leq p } \right ) \subset A_p
\]
and in addition $I_1 = (e^{xy}_1, e^{zy}_1) \subset A_1$. Let us stress that we
can take the expressions in (\ref{eqn!pfaffians}) for any given $j$ to write
the polynomials $e^{xy}_i$, $e^{zy}_i$ and $e^{y}_{ij}$.

\item For $0 \leq p \leq n-1$, define the homogeneous ideals
$J_p \subset A_p$ as follows:
\begin{enumerate}
\item $J_0 = (x_1,z_1) \subset A_0$;
\item $J_p = (x_{p+1},z_{p+1}, y_1, \dots ,y_p) \subset A_p$, for $p\geq 1$.
\end{enumerate}

\item Finally, notice that $I_p \subset J_p \subset A_p$.
Set $R_p = A_p/I_p$ and consider $J_p$ as an ideal of $R_p$.

\end{enumerate}

\begin {defn}   \label {dfn!binomfaffiansideal}
We call the ideal $I_n$ of $A_n$ the generic
   $\binom {n} {2}$ Pfaffians ideal.
\end {defn}

\end{notation}

The main aim of this section is to prove by induction on $p =0,1, \dots ,n$
that $R_p=A_p/I_p$ is a Gorenstein graded ring whose dimension is equal to
$\dim A_0-1$ (hence $I_p$ has codimension $p+1$ in $A_p$). Our strategy is to
establish inductively that $R_p$ is the result of serial unprojection of type
Kustin--Miller (\cite {PR} Definition~1.2). Our main algebraic result is the
following theorem, which we will prove in Subsection~\ref{sub!pfofalgebrathm}.

\begin {theorem}  \label {thm!mainalgebrathm}
Let all the notation be as above.
\begin{enumerate}\renewcommand{\labelenumi}{(\alph{enumi})}
\item  For $p =1, \dots , n$,  $R_p$ is the unprojection ring of type
Kustin--Miller of the pair $J_{p-1} \subset R_{p-1}$.
\item For $p =0,1, \dots , n$, $R_p$ is a normal Gorenstein graded 
integral domain, of dimension
equal to $\dim R_0$ (which is equal to  $2n+2^n$ since $\dim \Z =1$).

\item There are natural inclusions
\[
   R_0 \subset R_1 \subset \dots \subset R_p \subset K(R_0)
\]
where $K(R_0)$ is the field of fractions of $R_0$, all except the last induced
by the chain of  inclusions $A_0 \subset A_1 \subset \dots \subset A_p$.

\item For $p =0, \dots, n$ there exists a Zariski closed subset $F_p \subset
\Spec R_p$, with the codimension of $F_p$ in $\Spec R_p$ at least two such that
the open subscheme $\Spec R_p \setminus F_p$ is naturally isomorphic with an
open subscheme of $\Spec R_0$.

\item  For $p =1, \dots, n$,  $x_p, z_p$ is a regular sequence of $R_p$.

\item  For $p =0,1, \dots, n$ and $1 \leq i < j \leq n$, $x_i,x_j$ is a regular
sequence of $R_p$.
\end{enumerate}
\end{theorem}

\begin {rem}
The most important conclusions of Theorem~\ref{thm!mainalgebrathm}  are
(a) and (b). However, for the purposes of the inductive step we need all six
statements.
\end{rem}

\begin{rem}
 For $1 \leq t \leq n$ the inclusion $R_{t-1}\subset R_t$ of (c) of
Theorem~\ref{thm!mainalgebrathm} is given by $R_t = R_{t-1}[s_t]$ where $s_t\in
K(R_0)$ is the rational function given by
\[
s_t =  \frac {x_ty_t-e^{xy}_t } {x_t} =   \frac {z_ty_t -e^{zy}_t} {z_t} 
\]
and 
\[
    s_t = \frac {y_iy_t -e^{y}_{it}} {y_i} 
\]
for $1 \leq i \leq t-1$.
\end{rem}

\begin {rem}  \label {rem!commutativityofunprojections}
Fix $2 \leq p \leq n$.
Inside $\Spec R_0$ we have the $p$ codimension one subschemes
$V(x_i,z_i)$  for $1 \leq i \leq p$.
We can interpret
Theorem~\ref{thm!mainalgebrathm} and  Corollary~\ref{cor!aboutidentf} below
as  saying that the order we perform  the unprojections of the
subschemes is irrelevant. An interesting open question is to find general
conditions that will guarantee this kind of commutativity of unprojections.
\end{rem}

\subsection {Proof of Theorem~\ref{thm!mainalgebrathm}}  \label {sub!pfofalgebrathm}

We begin the proof of Theorem~\ref{thm!mainalgebrathm}.
We will need the following proposition.

\begin {prop}  \label{prop!forxpzpregseq }
Fix $1 \leq p \leq n$. Assume  $R_p$ is Cohen--Macaulay with $\dim R_p = \dim R_0$.
Then $x_p, z_p$ is a regular sequence for $R_p$.
\end {prop}

\begin {pf}
   Denote by $T \subset A_p$ is the ideal of $A_p$  generated by all $r_{d_1 \cdots d_n}$,
for  $(d_1, \dots ,d_n) \in \{ 0,1 \}^n $  with  $ (d_1, \dots ,d_n) \not= (0,0, \dots ,0)$ and
$ (d_1, \dots ,d_n) \not= (1,1, \dots ,1)$.

    Using  the assumptions about the dimension and the Cohen--Macaulayness of $R_p$,
to prove the proposition it  is enough to show that
\[
     \dim R_p/(x_p,z_p) \leq  \dim R_0 -2 =  \dim A_0 - 3
\]
and for that it is enough to show that
\[
     \dim A_p /(I_p + (x_p,z_p,y_p) + T)  \leq (\dim A_0 - 3) -(2^n-2)-1 = 2n-1  .
\]
We denote by $\commandeta (e^{xy}_i)$ the result of substituting
to  $ e^{xy}_i$ zero for $x_p, z_p, y_p$ and  $r_{d_1 \cdots d_n}$ for  $(d_1, \dots ,d_n) \in \{ 0,1 \}^n $
with  $ (d_1, \dots ,d_n) \not= (0,0,\dots ,0)$ and   $ (d_1, \dots ,d_n) \not= (1,1, \dots ,1)$,
and similarly for $\commandeta (e^{zy}_i)$ and $\commandeta (e^{y}_{ij})$. 

For $1 \leq i < p$ we have
\[
    \commandeta (e^{xy}_i)  = x_iy_i, \quad \quad  \commandeta (e^{zy}_i)  = z_iy_i,
\]
for $1 \leq i < j  < p$ we have
\[
    \commandeta (e^y_{ij})  = y_iy_j,
\]
and for $ 1 \leq i <  p$ we have
\begin {eqnarray*}
  & &    \commandeta (e^y_{ip}) = r_{00 \cdots 0} r_{11 \cdots 1}
               \Biggl[ \prod_{t=1}^{i-1} x_tz_t   \Biggr]       \Biggl[  \prod_{t=i+1}^{p-1} x_tz_t \Biggr]
                    \Biggl[ \prod_{t=p+1}^{n} x_tz_t   \Biggr]  \\
  & &    \commandeta (e^{xy}_p) = r_{11 \cdots 1}
               \Biggl[ \prod_{t=1}^{p-1}z_t   \Biggr]    \Biggl[   \prod_{t=p+1}^nz_t  \Biggr]  \\
  &  & \commandeta (e^{zy}_p) = -r_{00 \cdots 0}
                 \Biggl[  \prod_{t=1}^{p-1}x_t  \Biggr]   \Biggl[ \prod_{t=p+1}^nx_t \Biggr] .  \\
\end {eqnarray*}
For the proof of the first equality, substitute $x_p=z_p=0$ to 
\[
   e^{xy}_i =  x_iy_i + x_p  \dd{^2{Q'}}{z_i\partial x_p} + z_p  \dd{^2{Q'}}{z_i\partial z_p}
\]
where  $Q' =  r_{00 \cdots 0} \; x_1 \dots x_n + r_{11 \cdots 1} \; z_1 \dots z_n$.
The proof of the other equalities is similar.
Since
\[
   I_p = (e^{xy}_1, \dots , e^{xy}_p, e^{zy}_1, \dots , e^{zy}_p, e^{y}_{ij}) \subset A_p
\]
(with indices $1 \leq i < j \leq p$), using the vanishing of
$ \commandeta (e^y_{ij})=y_iy_j$, for $1 \leq i < j  < p$, we get two cases.

Case 1. All $y_i$ are zero, for $1 \leq i \leq p-1$. Then,  the vanishing of $\commandeta(e^{xy}_p)$ and
$\commandeta(e^{zy}_p)$  imply that two more variables vanish, so we get the desired codimension.

Case 2. There exists unique nonzero $y_a$, with $1 \leq a \leq p-1$. Using the vanishing
of   $\commandeta (e^{xy}_a)$ and $\commandeta (e^{zy}_a)$ we get
the vanishing of both $x_a$ and $z_a$, and using the vanishing of $ \commandeta (e^y_{ap})$ we
get that at least one more variable should vanish, so we again reach the desired codimension
which finishes the proof of the proposition.
\QED \medskip
\end {pf}

We now start the induction for the proof of Theorem~\ref{thm!mainalgebrathm}.

\begin {lemma}  \label{lemma!thmforpzero}
     Theorem~\ref{thm!mainalgebrathm} is true for $p =0$.
\end {lemma}

\begin {pf}
The ring $A_0$ is a Gorenstein normal integral domain, since it is
a finitely generated polynomial $\Z$-algebra. Since $Q \in A_0$ is an irreducible
polynomial, it follows that $R_0 = A_0 / (Q)$  is a Gorenstein integral domain.
Therefore, to prove that, for fixed $1 \leq i < j \leq n$, 
$x_i,x_j$ is a regular sequence of $R_0$ it is enough to 
show $\dim R_0/(x_i,x_j) = \dim R_0/(x_i,x_j) -2$, which
follows from the fact that  $Q$ does not vanish if we substitute 
to it $x_i = x_j =0$.
We will prove the normality of 
of $R_0$ by applying \cite {BV} Lemma~16.24.

Fix $1 \leq i \leq n$. The ring $R_0/(x_i)$ is a reduced ring, 
since the polynomial obtained from $Q$ by substituting zero for $x_i$
has no multiple factors.

Denote, for $1 \leq i \leq n$, by $T_i \subset R_0$ the multiplicatively closed subset
\[
   T_i = \{ 1, x_i, x_i^2, x_i^3,  \dots  \} \subset R_0.
\]
For notational convenience we also set $T_0 = \{ 1 \} \subset R_0$.
We will prove by induction on $t= 0,1, \dots ,n$ that 
the localization ring 
\[
   B_t = T_{n-t}^{-1} T_{n-t-1}^{-1} \dots  T_1^{-1} T_0^{-1} R_0
\]
is a normal integral domain.  Since a localization of an integral domain
is an integral domain, we only need to prove the normality of $B_t$.

Assume first that $t=0$.  By the form of $Q$, $B_0$ is isomorphic to a
localisation of the polynomial subring
\[
    \Z [x_1, \dots ,x_n, z_1, \dots ,z_n, r_{d_1 \dots d_n}] \subset A_0
\]
where $(d_1,d_2, \dots ,d_n) \in \{ 0,1 \}^{n}$  and 
$(d_1, \dots ,d_n) \not= (0,0, \dots ,0)$. Since 
the localisation of a normal ring is again normal, we get that
$B_0$ is a normal domain.

Assume now that for some $t$ with  $0 \leq t \leq n-1$ we have
that $B_t$ is normal. By \cite {BV} Lemma~16.24, to prove that
the domain $B_{t+1}$ is normal it is enough to prove that 
$B_{t+1}/(x_{n-t})$ is reduced. Since localisation commutes with taking
quotients and the localisation of a reduced ring is again reduced we have
that $B_{t+1}/(x_{n-t})$ is reduced as a localisation of the
already proven reduced $R_0/(x_{n-t})$. This finishes the induction 
and  hence the case $p=0$ of Theorem~\ref{thm!mainalgebrathm} follows.
\QED \medskip
\end {pf}

\begin {lemma}  \label{lemma!thmforpone}
   Theorem~\ref{thm!mainalgebrathm} is true for $p =1$.
\end {lemma}

\begin {pf}
 Using \cite{P2} Section~4, we get that $R_1$ is the unprojection of
type Kustin--Miller of the pair $J_0 \subset R_0$, hence by
the definitions of unprojection (\cite {PR} Section~1)
we have that $R_0$ is contained in a natural way in $R_1$ and $R_1$ has the
same dimension as $R_0$. Moreover, by  \cite{PR} Theorem~1.5 $R_1$ is
Gorenstein, and by \cite{PR} Remark~1.5, $R_1$ is a domain contained
in a natural way in the field of fractions $K(R_0)$ of $R_0$.

Proposition~\ref{prop!forxpzpregseq } says
that $x_1,z_1$ is a regular sequence of $R_1$, hence by setting
\[
       F_1 = V(x_1,z_1) \subset \Spec R_1
\]
we get that  $F_1$ has codimension two in $\Spec R_1$,
and by the construction of unprojection
$\Spec R_1 \setminus F_1$ is isomorphic in a natural way
(induced by the inclusion $R_0 \subset R_1$) to the open
subscheme   $\Spec R_0 \setminus V(x_1,z_1)$ of $\Spec R_0$.
Using Serre's normality criterion (\cite {BH} Theorem~2.2.22),
we get that the integral domain $R_1$ is normal.

We now prove that if $1 \leq i < j \leq n$ then
$x_i,x_j$ is a regular sequence of $R_1$.
If this was not true, using that we proved that $R_1$ is Gorenstein,
we would have that $V(x_i,x_j) \subset \Spec R_1$ would have
codimension at most one in $\Spec R_1$.
Using  the natural isomorphism of $\Spec R_1 \setminus F_1$
with the open subscheme of $\Spec R_0 \setminus V(x_1,z_1)$
and  that we proved that $F_1 \subset  \Spec R_1$ has codimension two in
$\Spec R_1$, it follows that $x_i,x_j$ is not
a regular sequence for $R_0$, contradicting Lemma~\ref{lemma!thmforpzero}.
 This finishes the proof of  Theorem~\ref{thm!mainalgebrathm}
for $p =1$.
\QED \medskip
\end {pf}

We now do the inductive step in the proof of  Theorem~\ref{thm!mainalgebrathm}.
We fix  $q$ with  $1 \leq q \leq n-1$. We assume that Theorem~\ref{thm!mainalgebrathm}
is true for all values $p$ with $0 \leq p \leq q$ and we will show
that it is true also for the case $p = q+1$.

For the rest of the proof,  given $0\leq t\leq n-1$, we will
denote  by $L_t \subset R_0$ the ideal
$L_t = (x_{t+1},z_{t+1}) \subset R_0$, and by $i_t \colon J_t \to R_t$  
and $i_{1,t} \colon L_t \to R_0$ the natural inclusion homomorphisms.

\begin {lemma}   \label{lem!phiexists}
There exists a unique homomorphism of Abelian groups
\[
   \phi_{q} \colon  \Hom_{R_0}(L_{q}, R_0) \to
            \Hom_{R_{q}}(J_{q}, R_{q})
\]
such that
\[
     \phi_{q}(f)  (x_{q+1}) = f(x_{q+1})
\]
for all  $f \in \Hom_{R_0}(L_{q}, R_0)$.
\end {lemma}

\begin {pf}
Recall that if $L \subset R$ is an ideal of a commutative ring $R$ and
$w \in L$ a nonzero divisor of $R$, then $\Hom_R(L, R)$ is isomorphic
to the ideal $\{a \in R \colon aL \subset (w) \}$ by the map
$ f \mapsto f(w)$ (cf.~\cite {PR} Remark 1.3.3). In particular,
$f$ is uniquely specified by the value $f(w)$.

Accordingly, since by the inductive hypothesis both $R_0$ and $R_{q}$
are integral domains with $R_0 \subset R_{q}$, it is enough to show that
\[
  \{a \in R_0 \colon  aL_{q} \subset R_0x_{q+1} \}  \subset
     \{b \in R_{q} \colon  bJ_{q} \subset R_{q}x_{q+1} \}
\]
Suppose $a \in R_0$ and $aL_{q} \subset R_0x_{q+1}$. In particular,
$az_{q+1} \in R_0x_{q+1}$. Obviously $a \in  R_{q}$, so we
need to show that $ay_i \subset R_{q}x_{q+1}$ for all $1 \leq i \leq q$.

Fix $1 \leq i \leq q$. Using the equation $e^{xy}_{i}$, which is zero
at $R_{q}$, we get
\begin {equation}  \label {eqn!forexistanceofmap}
        x_iay_i = - (ax_{q+1}Q_{i,{q+1}}^{zx}+az_{q+1}Q_{i,{q+1}}^{zz})
 \in R_{q} x_{q+1}
\end {equation}
By the inductive hypothesis, $x_i, x_{q+1}$ is a regular
sequence for $R_{q}$. As a consequence,  (\ref {eqn!forexistanceofmap})
implies that $ay_i \in R_{q} x_{q+1}$. Hence we get the existence of
the map $\phi_{q}$.  $R_{q}$ is a domain which implies that
$x_{q+1}$ is a regular element of $R_{q}$. The
uniqueness follows by the fact that an
element of $\Hom_{R_{q}}(J_{q}, R_{q})$ is uniquely specified
by its value at  $x_{q+1}$.
\QED \medskip
\end {pf}

Notice that clearly  $\phi_{q}(i_{1,q}) = i_{q}$.

\begin {prop}  \label {prop!existcnceofb}
Assume $f \in \Hom_{R_{q}}(J_{q}, R_{q})$. There
exists $b \in R_{q}$ such that the homomorphism
$f - bi_{q}$ maps $x_{q+1}$ and  $z_{q+1}$
inside $R_0 \subset R_{q}$.
\end {prop}

\begin {pf}
We prove by induction that for every $t =0, \dots ,q$ there exists
$b_t \in R_{q}$ such that $f - b_t i_{q}$ maps the elements
$x_{q+1}, z_{q+1}$ and  $y_j$, for $1 \leq j \leq q-t$, inside
$R_{q-t} \subset R_{q}$.

The result is trivially true when $t=0$. Assume $0 \leq t \leq q-1$
and that there exists $b_{t} \in R_{q}$ such that $f - b_ti_{q}$
maps the elements $x_{q+1}, z_{q+1}$ and  $y_j$, for $1 \leq j \leq q-t$,
inside $R_{q-t}$. Since, by construction,
$R_{q-t} = R_{q-(t+1)}[y_{q-t}]$ (as algebras),  there exist
$a \in  R_{q-(t+1)}$ and $c \in R_{q}$ with
\[
        (f - b_ti_{q}) (y_{q-t}) = a + cy_{q-t}.
\]
Set $g = f - (b_t+c)i_{q}$ We claim that $g$ maps  the elements
$x_{q+1}, z_{q+1}$ and  $y_j$, for $1 \leq j \leq q-(t+1)$,
inside $ R_{q-(t+1)}$. Indeed, if $u$ is in the ideal of $R_{q-(t+1)}$
generated by $x_{q+1}, z_{q+1}$ and  $y_j$, for $1 \leq j \leq q-(t+1)$,
we have
\[
    y_{q-t} g(u) =  u g(y_{q-t}) \in  R_{q-(t+1)}.
\]
Since by the inductive hypothesis of Theorem~\ref{thm!mainalgebrathm}
we have  normality of  $R_{q-(t+1)}$ and  that
$R_{q-t}$ is an unprojection of $R_{q-(t+1)}$, using
\cite {PR} Remark~1.3.4 (cf.~\cite{P1} Lemma~2.1.7) we get that
$g(u) \in R_{q-(t+1)}$, which finishes the proof of Proposition~\ref{prop!existcnceofb}.
\QED \medskip
\end {pf}

\begin {cor}  \label {cor!aboutidentf}
Fix $s \in  \Hom_{R_0}(L_{q}, R_0)$ such that $s$ together with $i_{1,q}$
ge\-ne\-ra\-te the $R_0$-module $\Hom_{R_0}(L_{q}, R_0)$. Then $\phi_{q}(s)$
together with $i_{q}$ generate the $R_{q}$-module
$\Hom_{R_{q}}(J_{q}, R_{q})$.
\end {cor}

\begin {pf}  Assume $f \in \Hom_{R_{q}}(J_{q}, R_{q})$.
Using Proposition~\ref{prop!existcnceofb}, there exists $b \in R_{p_{0}}$ such that,
if we set $g = f - bi_{q}$, we have $g(x_{q+1}) \in R_0$ and $g(z_{q+1})  \in R_0$.
Therefore, there exists $g_1 \in \Hom_{R_0}(L_{q}, R_0)$, with
$g_1(x_{q+1})=g(x_{q+1})$ and $g_1(z_{q+1}) = g(z_{q+1})$. By the
assumptions there exists  $c_1,c_2 \in R_0$ with $g_1 = c_1s+c_2 i_{1,q}$.
As a consequence,
\[
    f(x_{q+1})  = c_1 \bigr[ \phi_{q}({s}) (x_{q+1}) \bigl]  + (b+c_1) i_{q}(x_{q+1})
\]
which implies that
\[
     f =   c_1\phi_{q}({s})+ (b+c_1) i_{q}
\]
and the result follows.
\QED \medskip
\end {pf}

\begin {prop}  \label {prop!rpoplus1isunprojection}
The ring $R_{q+1}$ is isomorphic to the unprojection ring of the pair $J_{q} \subset R_{q}$.
\end {prop}

\begin {pf}
To simplify the notation of the proof we set, for  $a,b \in \{ x,z \}, \; Q^{ab} = Q^{ab}_{i,q+1}$.
Using \cite {P2} Section~4, $\Hom_{R_0}(L_{q}, R_0)$
is generated as $R_0$-module by the inclusion map  $i_{1,q}$ together
with the homomorphism  $t \colon L_q  \to R_0$ such that
\begin {eqnarray*}
  t(x_{q+1}) =  - (x_iQ^{xz} +  z_iQ^{zz}), \quad
   & & t(z_{q+1}) =    x_iQ^{xx} +  z_iQ^{zx}
\end {eqnarray*}
for all $1 \leq i \leq q$.
Notice that  these equations correspond exactly to $e^{xy}_{q+1}$ and $e^{zy}_{q+1}$.
Using Corollary~\ref {cor!aboutidentf}, $i_{q}$ together with $\; \phi_{q}(t)$ generate
the $R_{q}$-module  $\Hom_{R_{q}}(J_{q}, R_{q})$, so $\; \phi_{q}(t)$ can be used
to define the unprojection ring.

Fix $1 \leq i \leq q$.  We have inside $R_{q}$
\[
   x_{q+1}  \bigr[ \phi_{q}(t)(y_i) \bigl] = y_i \bigr[   \phi_q(t)(x_{q+1})\bigl] 
   = - y_i (x_iQ^{xz} +  z_iQ^{zz})
\]
Using the relations $e^{xy}_{i}=0$ and $e^{zy}_{i}=0$ which hold in $R_{q}$ (since
$1 \leq i \leq q$) we get
\begin {eqnarray*}
& &   - y_i(x_iQ^{xz} +  z_iQ^{zz}) =
             - Q^{xz} (x_iy_i) -  Q^{zz} (z_iy_i) = \\
  &  &   Q^{xz}  (x_{q+1}Q^{zx}+z_{q+1}Q^{zz})  - Q^{zz} (x_{q+1}Q^{xx}+z_{q+1}Q^{xz}) = \\
  &  &    x_{q+1} ( Q^{xz}Q^{zx}- Q^{xx}Q^{zz}).
\end {eqnarray*}
hence,  since $R_{q}$ is a domain,
\[
  \phi_{q}(t)(y_i) - (  Q^{xz}Q^{zx}- Q^{xx}Q^{zz})  = 0
\]
which corresponds exactly to $e^{y}_{iq}$.
As a consequence,  Proposition~\ref{prop!rpoplus1isunprojection} follows.
\QED \medskip
\end {pf}

\begin {prop}  \label {prop!rpoplusisgorenstein}
The ring $R_{q+1}$ is a Gorenstein integral domain, of
dimension equal to $\dim R_0$, containing $R_{q}$ in a
natural way and contained in the field of fractions
$K(R_0)$.
\end {prop}

\begin {pf}
 Using Proposition~\ref{prop!rpoplus1isunprojection}, we get
by  the definitions of unprojection (\cite {PR} Section~1)
that $R_{q}$ is contained in a natural way in $R_{q+1}$ and $R_{q+1}$ has the
same dimension as $R_{q}$, hence by the inductive hypothesis same
dimension as $R_0$.  Moreover, by  \cite{PR} Theorem~1.5 and
the inductive hypotheses for $R_{q}$ we get that  $R_{q+1}$ is
Gorenstein, and by \cite{PR} Remark~1.5 that it is also a domain contained
in a natural way in the field of fractions $K(R_{q})$ of $R_{q}$.
Since by the inductive hypothesis $K(R_{q}) = K(R_0)$,
Proposition~\ref{prop!rpoplusisgorenstein} follows.
\QED \medskip
\end {pf}

\begin {prop}  \label {prop!rpoplusoneregseq}
The sequence  $x_{q+1},z_{q+1}$ is a regular sequence
for $R_{q+1}$.
\end {prop}

\begin {pf}
By Proposition~\ref{prop!rpoplusisgorenstein}, $R_{q+1}$ is a Gorenstein
integral domain, of dimension equal to the $\dim R_0$. As a consequence, the
result follows by using  Proposition~\ref{prop!forxpzpregseq }.
\QED \medskip
\end {pf}

\begin {prop}  \label {prop!rpoplus1zariskiopen}
There exists a Zariski closed subset
$F_{q+1} \subset \Spec R_{q+1}$, with the codimension of $F_{q+1}$ in $\Spec R_{q+1}$
at least two such that the open subscheme $\Spec R_{q+1} \setminus F_{q+1}$
is naturally isomorphic with an open subscheme of $\Spec R_0$.
\end {prop}

\begin {pf}
By the construction of unprojection,
\[
  \Spec R_{q+1} \setminus V(x_{q+1}, z_{q+1}, y_1, \dots ,y_{q})
\]
is naturally isomorphic to an open subset of  $\Spec R_{q}$.
Using the inductive hypothesis and  Proposition~\ref{prop!rpoplusoneregseq}
the result follows.
\QED \medskip
\end {pf}

\begin {prop}  \label {prop!rpoplus11normality}
The ring $R_{q+1}$ is a normal domain.
\end {prop}

\begin {pf}
By Proposition~\ref{prop!rpoplusisgorenstein}, $R_{q+1}$  is a
Gorenstein integral domain. The result follows by combining the normality of
$R_0$ (Lemma~\ref{lemma!thmforpone}), Proposition~\ref{prop!rpoplus1zariskiopen}
and  Serre's normality criterion (\cite {BH} Theorem~2.2.22).
\QED \medskip
\end {pf}

\begin {prop}  \label {prop!rpoplus1refseqnumber2}
If $1 \leq i < j \leq n$ then
$x_i,x_j$ is a regular sequence of $R_{q+1}$.
\end {prop}

\begin {pf}
If this was not true, using that $R_{q+1}$ is Gorenstein
(Proposition~\ref{prop!rpoplusisgorenstein}) we would have that $V(x_i,x_j) \subset \Spec R_{q+1}$
would have codimension at most one in $\Spec R_{q+1}$.
Using  Proposition~\ref{prop!rpoplus1zariskiopen},
we get that $x_i,x_j$ is not
a regular sequence for $R_0$, contradicting Lemma~\ref{lemma!thmforpzero}.
\QED \medskip
\end {pf}

We have now finished the proof of the inductive step, hence
the proof of Theorem~\ref{thm!mainalgebrathm}.

\subsection {Generic perfection of $R_p$} \label {sub!genericperfection}

We fix $n \geq 2$ and $0 \leq p \leq n$. We will prove that the $A_p$-module $R_p$
is a generically perfect $A_p$-module. Recall (\cite {BV} Section~3.A) that this means
that $R_p$ is a perfect $A_p$-module and also faithfully flat as $\Z$-module.

A useful consequence of the generic perfection of $R_p$ is that whenever we substitute  
the variables of the ideal $I_p$ with elements of an arbitrary Noetherian ring 
we get, under mild conditions,  good induced properties of the resulting ideal 
(cf. \cite{BV} Section~3 for precise statements).  We will use the generic perfection 
of $R_p$ in Corollary~\ref{cor!gorensteinessoverk}, 
Remark~\ref{rem!aboutcomplexes} and Proposition~ \ref{prop!theresolution}.

\begin {rem} Recall (\cite {BV} Section~16.B) that if $A$ is Noetherian ring
and $M$ a finitely generated $A$-module, the grade of $M$ is defined to be
the maximal length of an $A$-regular sequence contained in the annihilator
ideal $\Ann M$ of $M$. If in addition $A$ is graded Cohen--Macaulay and $M$ is a
graded module, we have that the grade of $M$ is equal to $\dim R - \dim R/\Ann M$
(cf.~\cite {BH} Corollary 2.1.4). As a consequence, using Theorem~\ref{thm!mainalgebrathm},
$R_p$ has grade as $A_p$-module  equal to $p+1$.
\end {rem}

\begin {prop}  \label {prop!genericperfection}
The  $A_p$-module $R_p$ is generically perfect of grade  $p+1$.
\end {prop}

\begin {pf}
Using \cite{BV} Proposition~3.2 it is enough to prove that  $R_p$ is a perfect
$A_p$-module, and for every prime integer $p$ the $A_p \otimes_{\Z} \Z/p$-module
$R_p \otimes_{\Z} \Z/p$  is perfect.

Using \cite{BV} Proposition~16.19, the perfection of $R_p$ as $A_p$-module
follows from the Gorensteiness of $R_p$ (Theorem~\ref{thm!mainalgebrathm})
together with the fact that every finitely generated $A_p$-module has
finite projective dimension (cf.~\cite {BV} p.~35).

Fix an integer prime $p$. It is clear that all the arguments we used to
prove Theorem~\ref{thm!mainalgebrathm} work also if we replace $\Z$ by
$\Z/p$. As a consequence, we can argue as in the case of $R_p$ to get that
the $A_p \otimes_{\Z} \Z/p$-module $R_p \otimes_{\Z} \Z/p$  is perfect,
which finishes the proof of the proposition.
\QED \medskip
\end {pf}

\begin {cor}  \label {cor!gorensteinessoverk}
Let $k$ be an arbitrary field. The $A_p \otimes k$-module $R_p \otimes_{\Z} k$  is
perfect, of grade equal to $p+1$. Moreover the $k$-algebra $R_p \otimes_{\Z} k$
is Gorenstein.
\end {cor}

\begin {pf} It follows immediately by combining the Gorensteiness of $R_p$
(Theorem~\ref{thm!mainalgebrathm}) and  the generic perfection of $R_p$
(Proposition~\ref{prop!genericperfection}) with \cite {BV} Theorems 3.3 and 3.6.
\QED \medskip
\end {pf}

\begin {rem}  \label{rem!aboutcomplexes}
Using the construction of unprojection in \cite {KM} which is based on
resolution complexes, together with the fact that $J_p$ has Koszul complex as 
minimal resolution over $A_p$ (since it is generated by a regular sequence), we 
can inductively build the minimal graded resolution of $R_p$ over $A_p$. Using \cite {BV}
Theorem 3.3 this will give us the minimal graded resolution of $R_p\otimes_{\Z} k$ over
$A_p \otimes_{\Z} k$, where $k$ is an arbitrary field. We will use this remark in
Proposition~\ref{prop!theresolution}.
\end {rem}

In the following we fix an arbitrary field $k$, and set
\[
      A_p^k = A_p \otimes_{\Z} k,  \quad \quad  R_p^k = R_p \otimes_{\Z} k.
\]

\begin {lemma}  \label {lemma!regseqafterfield}
The length $n+2^n-1$ sequence  
$ \; z_1, \dots , z_n, r_{d_1 \cdots d_n}  \; $ 
with indices  $(d_1, \dots ,d_n) \in \{ 0,1 \}^{n}$  and 
$(d_1, \dots ,d_n) \not= (0,0, \dots ,0)$ is regular for
$R_p^k$ with respect to any ordering of it.
\end {lemma}

\begin {pf}
Denote by $T \subset R_p^k$ the ideal of $R_p^k$  generated by the sequence.
Since by Corollary~\ref{cor!gorensteinessoverk} $R_p^k$
is Gorenstein, hence Cohen--Macaulay, it is enough to prove that
\[
     \dim R_p^k/T =  \dim R_p^k - (n+2^n-1)
\]

Denote by $T_1$ the monomial ideal of $k[x_1, \dots ,x_n, y_1, \dots ,y_p, r_{00 \cdots 0}]$
generated by  $x_iy_i$ for $1 \leq i \leq p$, $y_iy_j$ for $1 \leq i < j \leq p$,
and $w_t$ for $1 \leq t \leq p$, with
\[
   w_t = r_{00 \cdots 0}  \Biggl[   \prod_{i=1}^{t-1} x_i    \Biggr]   \Biggl[ \prod_{i=t+1}^n x_i   \Biggr].
\]
Arguing as in the proof of  Proposition~\ref{prop!forxpzpregseq }
we get that
\[
   R_p^k/T  \iso k[x_1, \dots ,x_n, y_1, \dots ,y_p, r_{00 \cdots 0}]/T_1
\]
and, moreover, that the right hand side ring has the right dimension.
\QED \medskip
\end {pf}

\begin {prop}  \label {prop!regseqafterfield}
Denote by $W_1$ the $k$-vector subspace of $A^k_p$ spanned by all
$r_{d_1 \cdots d_n}$ with $(d_1, \dots ,d_n) \in \{ 0,1 \}^{n}$ and
$(d_1, \dots ,d_n) \not= (0,0, \dots 0)$.  Assume $1 \leq t \leq 2^n-1$
and that $l_1, \dots l_t$ are elements of $W_1$ which are $k$-linearly
independent. Then $l_1, \dots l_t$ is a regular sequence for $R_p^k$.
\end {prop}

\begin {pf}
Since by Corollary~\ref{cor!gorensteinessoverk} $R_p^k$
is Gorenstein, hence Cohen--Macaulay, it is enough to prove that
the dimension drops by $t$ when we divide $R_p^k$ by the ideal
generated by $l_1, \dots ,l_t$. This follows from
Lemma~\ref{lemma!regseqafterfield} after completing $l_i$ to
a basis of $W_1$ and dividing by the ideal generated by the basis
together with $z_1, \dots ,z_n$.
\QED \medskip
\end {pf}

\section {The numerical Campedelli surface construction } \label {sec!campedelli}

In this section we work over the field $k=\C$ of complex numbers.
We will use the algebra developed in Section~\ref {sec!algebrapart} in order to prove
the existence of numerical Campedelli surfaces with torsion group equal to $\Z/6$.

We define the polynomial ring
\[
    A_4^s = k [x_1, \dots ,x_4, z_1. \dots ,z_4 , y_1, \dots y_4]
\]
($s$ for specific), and we  assign degree $1$ to each variable
$x_i$ and $z_i$,  for $1 \leq i \leq 4$, and  degree $2$ to each variable
$y_i$, for  $1 \leq i \leq 4$.

Let $G$ be the cyclic group of   order 6 and denote by $g$ a generator of $G$.
We define  a linear action of $G$ on $A^s_4$ by
\begin {eqnarray*}
   & &  (gx_1, gx_2,  gx_3 )   = (-x_2, -x_3, -x_1),  \quad \quad \quad \quad gx_4 = -x_4, \\
   & &  (gz_1, gz_2,  gz_3 ) \;  = ( z_2,  z_3,  z_1), \quad \quad \quad \quad \quad \quad \; \; \; gz_4 = z_4,\\
   & &  (gy_1, gy_2,  gy_3 ) \;  = (-y_2, -y_3, -y_1), \quad \quad \quad \quad  \;  gy_4 = -y_4. 
\end {eqnarray*}

Consider the sixteen dimensional $k$-vector subspace $W_1$ of the degree four polynomials
of  $A_4^s$ spanned by the monomials $a_1a_2a_3a_4$ where
$a_i \in \{x_i, z_i \}$. It is easy to see that $W_1$ is $G$-invariant,
and that the vector subspace  $W_2 \subset W_1$ of $G$-invariant elements
is $4$-dimensional with $k$-basis $F_1,\dots ,F_4$, where
\begin {eqnarray*}
 & &  F_1 = x_1x_2x_3x_4, \; \; \; \; F_2 = (x_1x_2z_3 +  x_1z_2x_3 + z_1x_2x_3)z_4,  \\
 & &  F_3 =  (x_1z_2z_3 + z_1x_2z_3 +z_1z_2x_3)x_4,  \; \; \; \; F_4 = z_1z_2z_3z_4.
\end {eqnarray*}

Fix $(r_1, \dots,r_4) \in k^4$ nonzero. We set
\begin {equation}  \label{eqn!dfnofQs}
   Q^s = Q^s(r_t) =  \sum_{i=1}^{4}  r_i F_i \in A_4^s.
\end {equation}
The polynomial $Q^s$ is homogeneous of degree four. Similarly to item (\ref{item!105})
of Notation~\ref{not!bignotation}
for each $1\leq i< j\leq n$, let
\[
\begin{array}{c}
Q_{ij}^{s,xx}=\dd{^2Q^s}{x_i\partial x_j},\quad Q_{ij}^{s,xz}=
\dd{^2Q^s}{x_i\partial z_j},\quad  Q_{ij}^{s,zx}= \dd{^2Q^s}{z_i\partial x_j}\quad
\text{and}\quad Q_{ij}^{s,zz}= \dd{^2Q^s}{z_i\partial z_j}.
\end{array}
\]
We clearly have 
\[
  Q^s = x_ix_j Q_{ij}^{s,xx}+x_iz_j Q_{ij}^{s,xz}+ z_ix_jQ_{ij}^{s,zx}+
         z_iz_jQ_{ij}^{s,zz}.
\]
Consider, for $1 \leq i < j \leq 4$,
the $5 \times 5$ skew--symmetric matrix
\[
    M^s_{ij} =  \begin {pmatrix}
   0   & x_i & z_i & -x_j & -z_j \\
       &   0 &  y_j & Q^{s,zz}_{ij} & -Q^{s,zx}_{ij}\\
       &     &  0   & -Q^{s,xz}_{ij} & Q^{s,xx}_{ij} \\
       &  \text{-sym}   &       & 0   & -y_i \\
       &    &      &   &   0
    \end {pmatrix}
\]
with entries in $A^s_4$. We denote by $I^s_4$ the ideal
of $A^s_4$ generated by all the submaximal Pfaffians of $M_{ij}$
for all values $1 \leq i < j \leq 4$.

The analogue of  item (\ref{item!101}) of Notation~\ref{not!bignotation} is true,
and we  denote, for
$1 \leq i < j \leq 4$,  by $e^{sy}_{ij}$ the Pfaffian of
$M^s_{ij}$  involving $y_iy_j$ with coefficient $1$,  by  $e^{sxy}_i$ the Pfaffian of
$M^s_{1i}$ (or of $M^s_{12}$ if $i=1$) involving $x_iy_i$ with coefficient $1$,  and
by  $e^{szy}_i$ the Pfaffian of  $M^s_{1i}$ (or of $M^s_{12}$ if $i=1$) involving
$z_iy_i$  with coefficient $1$.  See (\ref{eqn!explicits}) below for the explicit formulas 
of $e^{sxy}_i, e^{szy}_i$ and $e^{sy}_{ij}$.

It is also clear that
\begin {equation}  \label{eqn!gensofI5s}
   I^s_4 = ( e^{sy}_{ij}, e^{sxy}_t, e^{szy}_t ) \subset A^s_4
\end {equation}
with indices $1 \leq i < j \leq 4$ and $1 \leq t \leq 4$.

We denote $R^s_4 = A^s_4/I^s_4$ the quotient ring (which, of course,
also depends on the choice of parameter values $(r_t)$).

\begin {prop}  \label{prop!aboutrs4}
a)  For any $(r_1, \dots ,r_4)  \in k^4$, $\dim R^s_4 \geq 7$.

b) Whenever $\dim R^s_4 = 7$, $R^s_4$ is a
Gorenstein ring and a perfect $A^s_4$-module.

c) There exist parameter values $(r_t)$ such that $\dim R^s_4 = 7$.

d) For general parameter values  $(r_1, \dots ,r_4)  \in k^4$
(in the sense of being outside a proper Zariski closed subset of $k^4$)
the ring  $R^s_4$  is Gorenstein with $\dim R^s_4 = 7$.
\end {prop}

\begin {pf}
Part a) follows by combining Theorem~\ref{thm!mainalgebrathm}, 
Proposition~\ref{prop!regseqafterfield} and
\cite {L} Theorem~4.3.12.  Part b) follows by combining Theorem~\ref{thm!mainalgebrathm},
Proposition~\ref{prop!regseqafterfield} and \cite {BV} Theorem~3.5.

For part c) we make specific choice of
parameters $r_1=1$ and $r_i = 0$ for $2 \leq i \leq 4$. We will
prove that  $\dim R^s_4 = 7$. For that, it is enough to prove
that  $\dim R^s_4/(z_1,z_2,z_3,z_4) = 3$.
Arguing as in the proof of Proposition~\ref{prop!forxpzpregseq }, it is
easy  to see  (compare also (\ref{eqn!explicits})), that
\[
       R^s_4/(z_1,z_2,z_3,z_4)  \iso k[x_1, \dots,x_4, y_1, \dots ,y_4] / T,
\]
where $T$ is the monomial ideal of $k[x_1, \dots,x_4, y_1, \dots ,y_4]$ generated by
the elements $x_iy_i$, for $1 \leq i \leq 4$,  
together with $y_iy_j$, for $1 \leq i < j \leq 4$, together with $w_t$, for $1 \leq t \leq 4$,
where 
\[
   w_t = \Biggl[   \prod_{i=1}^{t-1} x_i    \Biggr]   \Biggl[ \prod_{i=t+1}^4 x_i   \Biggr]
\]
and that we indeed have the right dimension.

Using semicontinuity of the fiber dimension  (cf.~\cite {Ei} Corollary~14.9)
and parts a) b) c),  we have that part d) follows, which finishes the proof
of  Proposition~\ref{prop!aboutrs4}.
\QED \medskip
\end {pf}

\begin {rem}
A different way of arguing for the proof of Proposition~\ref{prop!aboutrs4}
is to suitably modify the arguments used in the proof of
Theorem~\ref{thm!mainalgebrathm}. One should be able
to get this way the more precise result that  $R^s_4$  is Gorenstein
with $\dim R^s_4 =7$ whenever there exists $i$, with $1 \leq i \leq 4$,
with $r_i$ nonzero. We will not use that in the following.
\end {rem}

Denote by $\zeta \in k$ a fixed primitive $6$th root of unity. We consider
the following homogeneous elements $m^i_j \in A_4^s$ of degree $1$
\begin {eqnarray*}
    m^0_1 & = & z_1+z_2+z_3,   \quad     m^0_2 = z_4,  \quad  m^1_1 = x_1+\zeta^2x_2 + \zeta^4x_3,  \\
    m^2_1 & = & z_1+\zeta^4z_2+\zeta^2z_3, \quad  m^3_1  = x_1+x_2+x_3, \quad  m^3_2 = x_4,  \\
    m^4_1 & = &  z_1+\zeta^2z_2+\zeta^4z_3,  \quad  m^5_1 = x_1+\zeta^4x_2 + \zeta^2x_3.
\end {eqnarray*}
By construction, each $m^i_j$ is an eigenvector for the action of $g \in G$ with eigenvalue
equal to $\zeta^i$, that is
\[
      g m^i_j = \zeta^i m^i_j.
\]
We fix four more complex numbers $(r_5, \dots , r_8) \in k^4$,  and we define
four homogeneous elements $h_i = h_i (r_5, \dots ,r_8) \in A_4^s $, for
$1 \leq i \leq 4$, by
\begin {eqnarray}
    h_1 &  = & m^0_1, \quad      h_2 = m^0_2, \quad   h_3 = m^3_2+r_5m^3_1,   \label {eqn!dfnofhi}\\
    h_4 & = &  y_4 + r_6 (y_1+y_2+y_3)+ r_7 m^1_1m^2_1+   r_8 m^4_1m^5_1.  \notag
\end {eqnarray}
We have that, for $1 \leq i \leq 4$,  the element $h_i$ is an eigenvector for $g$
with eigenvalue equal to $1$ for $h_1$ and $h_2$ and eigenvalue equal to $-1$ for
$h_3$ and $h_4$.

We  denote by $T^s  \subset A_4^s$ the homogeneous ideal
\[
     T^s =  T^s(r_t) = I_4^s + (h_1, \dots ,h_4) \subset A_4^s
\]
Moreover, we denote by $A$ the polynomial subring
\[
     A =  k [x_1,x_2,x_3, z_1,z_2,y_1,y_2,y_3] \subset A_4^s
\]
with the weighting of the variables induced by that of $A_4^s$.
For fixed general parameter values  $(r_1, \dots ,r_8) \in k^8$,
the composition
\[
     A \to A_4^s   \to A_4^s / T^s
\]
is surjective (where the first map is the natural inclusion and the
second the natural projection), so we get an induced isomorphism
\begin {equation}   \label {eqn!isoofgradedrings}
      \frac{ A } {L} \iso    \frac{ A^s_4} {T^s}
\end {equation}
where $L \subset A$ is the kernel of the composition.

\begin {prop}   \label {prop!hiregularsequence}
a)  For any choice of parameter values $(r_1, \dots ,r_8) \in k^8$ we
have $ \dim A_4^s/T^s \geq 3$ and  whenever $\dim A_4^s/T^s = 3$
we have that $A_4^s/T^s$ is a Gorenstein ring, perfect as $A_4^s$-module.

b) There exist parameter values $(r_t)$ such that $ \dim A_4^s/T^s = 3$.

c)  For general parameter values  $(r_1, \dots ,r_8) \in k^8$,
(in the sense of being outside a proper Zariski closed subset of $k^8$)
$ \dim A_4^s/T^s = 3$ and  $A_4^s/T^s$  is a Gorenstein ring.

\end {prop}

\begin {pf}
Part a) follows immediately from Proposition~\ref{prop!aboutrs4}
by noticing that $A_4^s/T^s$ is isomorphic to $R_4^s/(h_1, \dots ,h_4)$.

For part b)  we fix the parameter values $r_1=r_4 =1$ and $r_j =0$, for
$2 \leq j \leq 8$ with $j \not= 4$.
By (\ref{eqn!isoofgradedrings}) $ A_4^s/T^s  \iso A/L$,
where $L$ is the ideal of  $A$ generated by 
\begin {eqnarray*}
 & & \{ x_1y_1,z_1y_1,x_2y_2,z_2y_2,x_3y_3,z_1y_3+z_2y_3,z_1^2z_2+z_1z_2^2, x_1x_2x_3, y_1y_2,  \\
& &  \; \; \; y_1y_3, x_2x_3z_1z_2+x_2x_3z_2^2,y_2y_3,x_1x_3z_1^2+x_1x_3z_1z_2, x_1x_2z_1z_2 \}.  \notag
\end{eqnarray*}

It is easily checked that each minimal associated prime of $L$  has codimension five
in  $A$, hence
\[
   \dim  A/L  = 3.
\]

Part c) is an immediate consequence of parts a) and b) arguing as in the proof
of Proposition~\ref{prop!aboutrs4}.
\QED \medskip
\end {pf}

\begin {prop}  \label {prop!theresolution}
 For general parameter values  $(r_1, \dots ,r_8) \in k^8$ the
minimal graded resolution of  $A/L$ as $A$-module is equal to
\begin {eqnarray}
   0 \to  A(-12) \to  A(-9)^8 \oplus A(-8)^6  \to  A(-8)^3 \oplus A(-7)^{24} \oplus A(-6)^8  \to
                        \notag  \\
  \to  A(-6)^8 \oplus A(-5)^{24} \oplus A(-4)^3  \to  A(-4)^6 \oplus  A(-3)^8  \to  A  \quad 
          \quad \quad  \label{eqn!resolnofT}
\end {eqnarray}
Moreover, the dualising module of  $A/L$ is equal to $(A/L)(1)$ and the Hilbert series of
 $A/L$ as graded $A$-module is equal to
\[
   \frac {t^4+2t^3+6t^2+2t+1} {(1-t)^3}  \in \Q(t).
\]
\end {prop}

\begin {pf}
   Using Remark~\ref{rem!aboutcomplexes}, we can easily calculate inductively the minimal
graded resolution of the generically perfect  (Proposition~\ref{prop!genericperfection})
module $R_4$ over $A_4$. Equation~(\ref{eqn!resolnofT}) follows by combining
Proposition~\ref{prop!hiregularsequence},  \cite {BV} Theorem~3.5, and the
easily observed fact that the minimal graded resolution of $R_4$ over $A_4$ remains homogeneous
and minimal.  The other conclusions of Proposition~\ref{prop!theresolution} follow
easily from (\ref{eqn!resolnofT}).
\QED \medskip
\end {pf}

\begin {defn}  \label{dfn!dfnofS}
 For general $(r_1, \dots ,r_8) \in k^8$  we denote by $S$ to be
the scheme
\[
     S = S(r_t) = \Proj   A_4^s/T^s \subset \P(1^8,2^4).
\]
\end {defn}

Our main aim is  prove that $S$ is an irreducible nonsingular surface with
invariants $p_g = 5, q = 0, K^2 = 12$ and trivial algebraic
fundamental group, which is an \'{e}tale   six to one cover of a numerical Campedelli
surface.

\begin {rem} By (\ref{eqn!isoofgradedrings}),
$S$ has an embedding as a nondegenerate subscheme
\begin {equation}  \label{eqn!SinsmallP}
    S \subset \P(1^5,2^3).
\end {equation}
\end {rem}

\begin {prop}  \label {prop!variouspropertiesofS}
a) The homogeneous coordinate ring of the embedding
$S \subset \P(1^8,2^4)$ is isomorphic to  $A_4^s/T^s$.

   b)  The scheme $S$ is a projective purely two dimensional  scheme over
$k$. Moreover,  $S$ is connected and $H^1 (S,\Oh_S(t)) =0 $ for all $t \in \Z$.

   c) The dualising sheaf $\om_S$ is isomorphic to $\Oh_S(1)$ as  $\Oh_S$-module.
\end {prop}

\begin {pf}   The graded ring  $A_4^s/T^s$ is Gorenstein (Proposition~ \ref{prop!hiregularsequence}),
hence  sa\-tu\-ra\-ted.
As a consequence, part a) follows.

Using part a) the homogeneous coordinate ring of the embedding
$S \subset \P(1^8,2^4)$ is Gorenstein, hence Cohen--Macaulay. It is then
well known (cf.~\cite{Do}, \cite{Ei} Ch.~18)
that  the conclusions of part b) follow.  It is also well-known that part c)
follows immediately from Proposition~\ref{prop!theresolution}.
\QED \medskip
\end {pf}

In the following we will also need the affine cone  $S^c \subset \A^{12} $ over
$S \subset \P(1^8,2^4)$, so we set
\[
      S^c = V(T^s) \subset \A^{12}.
\]
We denote by  $S^c_{cl}$ the set of closed points of $S^c$, and by  $S_{cl}$ the set of closed
points of $S$.  Since $k = \C$ is algebraically closed, we can identify $S^c_{cl}$ with the set
of points
\begin {equation}  \label{eqn!forP}
  P= (a^x_1, \dots ,a^x_4; \; \; a^z_1, \dots ,a^z_4; \; \; a^y_1, \dots ,a^y_4)  \in k^{12}
\end {equation}
such that $f(P) = 0$ for every $f \in  T^s$.

By definition, $S_{cl}$ is the quotient of $S^c_{cl} \setminus \{0 \}$ under the
group action
\[
   k^* \times (S^c_{cl} \setminus \{ 0 \})  \to (S^c_{cl} \setminus \{ 0 \})
\]
with
\[
    hP  =  (ha^x_1, \dots ,ha^x_4; \; \; ha^z_1, \dots ,ha^z_4; \; \;  h^2a^y_1, \dots ,h^2a^y_4)
\]
for $h \in k^*$ and $P \in S^c_{cl}$ as in (\ref{eqn!forP}).

Since by Proposition~\ref{prop!invarianceofTs} below the ideal $T^s \subset R_4^s$
is $G$-invariant, there is an induced $G$ action  $G \times R_4^s/T^s \to R_4^s/T^s$,
which  induces in a natural way two group actions: $G \times S^c_{cl}  \to  S^c_{cl}$ and
$G \times S_{cl}  \to  S_{cl}$.

Explicitly, for $P  \in S^c_{cl}$ as in (\ref{eqn!forP}) we have
\begin {equation} \label{eqn!forgP}
      gP  =    (-a^x_3,-a^x_1,-a^x_2,-a^x_4; \; \; a^z_3,a^z_1, a^z_2,  a^z_4; \; \;
                     -a^y_3,-a^y_2,-a^y_1,-a^y_4 )
\end{equation}

\begin {lemma}   \label{lemma!nopointswithzeroxi}
For general values of parameters $(r_1, \dots ,r_8) \in k^8$
there is no nonzero point $P \in S^c_{cl}$ (notation for $P$ as in  (\ref {eqn!forP}))
such that $a^x_i = a^z_i = 0$ for all $1 \leq i \leq 4$.
\end {lemma}

\begin {pf}
Indeed, if all $a^x_i =  a^z_i = 0$ we have by looking at $e^{sy}_{ij}$, 
 for $ 1 \leq i < j \leq 4$ (cf.~(\ref{eqn!explicits})), that at least three of the four
$a^y_t$ are zero, and then by looking at the polynomial $h_4$ we get that the
remaining $a^y_t$ is also zero, a contradiction to $P \not=0$.
\QED \medskip
\end {pf}

\begin {prop}   \label {prop!nonsingularity}
Consider $S \subset \P(1^5,2^3)$  as in (\ref{eqn!SinsmallP}).  Denote by
$S^c_1 \subset \A^8$ the affine cone over $S$. The scheme $S^c_1$ is smooth
outside the vertex of the cone.
\end {prop}

\begin {pf}
Unfortunately, we were only able to prove Proposition~\ref{prop!nonsingularity}
with the help of the computer algebra program Singular \cite{GPS01}. We took
a similar approach  as in \cite{R2} p.~18 and worked over the finite field of
$\Z/103$ after putting values for parameters $r_7=r_8=0$, in order to have everything
defined over $\Z$.
\QED \medskip
\end {pf}

\begin {rem} Using the birational character of unprojection
it is not hard  to specify inductively  (for general values of the
parameters $(r_t)$)   the singularities of the affine cone over the
$6$-fold $V(I_t^s)$ for $t =0, \dots ,4$, where $I_t^s$ are the precise analogues
of the ideals $I_t$ defined in Section~\ref{sec!algebrapart}.
With a little more effort, one can also specify inductively the singularities of the 
cone over the $3$-fold  $V(I_t^s+ (h_1,h_2,h_3))$.
Since $h_2 = z_4$ vanishes, the trick here is to start from the codimension two ideal
$(e^{sxy}_{4}, e^{szy}_{4})$
and then inductively unproject $V(x_1,z_1,y_4),V(x_2,z_2,y_4,y_1)$ and finally 
$V(x_3,z_1+z_2,y_4,y_1,y_2)$.
What, unfortunately, we were not able to do
was to find a way to deduce the nonsingularity (outside the vertex of the
affine cone) of the surface from the singularity calculations of the
$3$-fold.
\end {rem}

\begin {theorem} \label{theorem!smoothnessandirreducibility}
Fix general values of parameters $(r_1, \dots ,r_8) \in k^8$.
$S=S(r_t)$ is an irreducible minimal nonsingular surface of general type with
$p_g = 5, q = 0, K^2 = 12$  and canonical ring isomorphic
to $ A_4^s/T^s$.
\end {theorem}

\begin {pf}
By combining Propositions~\ref{lemma!nopointswithzeroxi}  and  \ref{prop!nonsingularity} 
we get that the scheme $S$ is smooth. Since  $S$ is also connected
(Proposition~ \ref{prop!variouspropertiesofS}), it follows that $S$ is an
irreducible nonsingular surface.

By Proposition~\ref{prop!variouspropertiesofS} the dualising sheaf
$\omega_S $ is isomorphic to $\Oh_S(1) $.  Using Lemma~\ref{lemma!nopointswithzeroxi}
$\Oh_S(1)$ is globally generated. As a consequence, 
\[
   \Oh_S(1)^{\otimes n} \iso  \Oh_S(n)
\]
for all $n \geq 1$, hence $ A_4^s/T^s$ is isomorphic
to the canonical ring of $S$. Therefore $\omega_S$ is ample which implies that
$S$ is minimal.

Since the irregularity $q$ of $S$ is zero (because by
Proposition~\ref{prop!variouspropertiesofS} $h^1(S, \Oh_S) = 0$ ),
the properties $p_g = 5, K^2 = 12$ follow by comparing the Hilbert series calculation
of Proposition~\ref{prop!theresolution}   with  \cite{R1} Example~3.5.
\QED \medskip
\end {pf}

\begin {rem}  We will  prove below that $S$ has trivial algebraic  fundamental group
(see the proof of Theorem~\ref{theorem!existenceCampedelli}).
\end {rem}

Our next aim is to prove that $S$ is an \'{e}tale six to one cover of a numerical Campedelli surface.

\begin {prop}   \label {prop!invarianceofTs}
     Assume $g_1 \in G$ and $u \in T^s$, then $g_1u \in  T^s$.
\end {prop}

\begin {pf}
Since $G = \Span{g}$, it is enough to check that $gu \in T^s$, where $u$ is
one of the generators  of $I_4^s \;$ appearing in (\ref{eqn!gensofI5s}).

It is easy to check (compare~(\ref{eqn!explicits})) that for $i \in \{1,2,3 \}$ we have
\[
    ge^{sxy}_{i} = e^{sxy}_{t}, \quad \quad  ge^{szy}_{i} = -e^{szy}_{t}, 
      \quad  \quad  ge^{sy}_{i4} =  e^{sy}_{t4},
\]
where $t \in \{1,2,3 \}$ is uniquely specified by $ t \equiv i+1 \mod 3$,
and also that 
\[
  ge^{sxy}_{4} = e^{sxy}_{4},  \quad  ge^{szy}_{4} = -e^{szy}_{4}, \quad 
    ge^{sy}_{12} =  e^{sy}_{23}, \quad  ge^{sy}_{13} =  e^{sy}_{12}, \quad ge^{sy}_{23} =  e^{sy}_{13}.
\]

A more conceptual proof can be given by arguing that due to the $G$-invariance of
$Q^s$, the action of $g$ interchanges (up to sign) the set of $Q_{ij}^{s,ab}$
(for $a,b \in \{x,z \}$), and use that
to argue that the action of $g$ interchanges (up to sign differences of whole columns
or rows) the set of the matrices $M_{ij}^s$.
\QED \medskip
\end {pf}

The proof of the following proposition will be given in Subsection~\ref {sub!proofofbasepointfreeness}.

\begin {prop}   \label {prop!basepointpointfreeness}
Fix general values of the parameters $(r_1, \dots ,r_8)  \in k^8$.
If $g_1 \in G$ with $g_1$ not the identity element and $u \in S_{cl}$
we have  $g_1u \not= u$.
In other words, the action of $G$ on $S_{cl}$ is basepoint free.
\end {prop}

The following is our main result about the existence of numerical Campedelli surfaces with
algebraic fundamental group equal to $\Z/6$.

\begin {theorem}   \label {theorem!existenceCampedelli}
For general $(r_1, \dots ,r_8) \in k^8$, the action of $G$ on $S$ is basepoint free.
As a consequence, the quotient surface $S/G$ is a smooth irreducible minimal complex surface of general
type with  $p_g = q =0$ and $K^2 = 2$ (i.e., a numerical Campedelli surface).
Moreover $S/G$ has both algebraic fundamental group and torsion group isomorphic to $\Z/6$.
\end {theorem}

\begin {pf}
Fix general $(r_1, \dots ,r_8) \in k^8$. By Proposition~\ref{prop!basepointpointfreeness},
the action of $G$ on $S$ is basepoint free, hence using 
Theorem~\ref{theorem!smoothnessandirreducibility}, $S/G$ is a smooth irreducible surface.
Denote by $\pi \colon S \to S/G$ the natural projection map. Since  $\pi$ is \'{e}tale 
$\pi^*(\omega_{S/G}) \iso \omega_{S}$ (cf.~\cite{MP}~p.~3), and since by the proof
of  Theorem~\ref{theorem!smoothnessandirreducibility} $\omega_{S}$ is ample
we have that $\omega_{S/G}$ is ample (cf.~\cite{Ha}, Exerc.~III.5.7), hence 
$S/G$ is a minimal surface of general type.

The invariants of $S/G$ follow from those of $S$ calculated in
Proposition~\ref{theorem!smoothnessandirreducibility}. Indeed, $\pi$ surjective
and $q(S) =0$ imply $q(S/G) =0$, and $\pi$   \'{e}tale six to one 
imply  $K_{S}^2 = 6 K_{S/G}^2$ and $\chi(S) = 6\chi(S/G)$.

To prove that the algebraic fundamental group of $S/G$ is equal to $G$ 
it is enough to show that  $\pi_1^{\mathrm{alg}}S = 0$.
Assume it is not, then the group  $\pi_1^{\mathrm{alg}}(S/G)$ has $6 |\pi_1^{\mathrm{alg}}S| \geq 12$ elements,
which contradicts that a Campedelli surface has algebraic fundamental
group consisting of at most $9$ elements (cf.~\cite{BPHV} Chap. VII.10). 

Since  the torsion group of $S/G$ is the largest abelian quotient of 
$\pi_1^{\mathrm{alg}}(S/G)$ (cf.~\cite {MP} p.~16), we get that 
$S/G$ has torsion group isomorphic to $\Z/6$ which
finishes the proof of Theorem~\ref{theorem!existenceCampedelli}.
\QED \medskip
\end {pf}

\begin {rem}   \label{rem!forgetaction}
If one is not interested in a group action, a bigger family
of surfaces of general type with $p_g = 5, q = 0, K^2 = 12$ 
can be obtained
by setting in (\ref{eqn!dfnofhi})  $h_1, \dots ,h_4$ to
be general homogeneous elements of $A_4^s$ of degrees, respectively, $1,1,1,2$, 
and $Q^s$ in (\ref{eqn!dfnofQs}) to be a general $k$-linear combination 
of the sixteen degree $4$ monomials $a_1a_2a_3a_4$, where $a \in \{x,z \}$.
\end {rem}

\subsection {The proof of Proposition~\ref {prop!basepointpointfreeness} } \label{sub!proofofbasepointfreeness}

For the proof of Proposition~\ref {prop!basepointpointfreeness}  we will need the following formulas.
\begin {eqnarray}   
&  & e^{sxy}_{1} =  x_2x_3z_4r_2+x_3x_4z_2r_3+x_2x_4z_3r_3+z_2z_3z_4r_4+x_1y_1  \notag \\
&  & e^{szy}_{1}  = -x_2x_3x_4r_1-x_3z_2z_4r_2-x_2z_3z_4r_2-x_4z_2z_3r_3+z_1y_1 \notag  \\
&  & e^{sxy}_{2} =  x_1x_3z_4r_2+x_3x_4z_1r_3+x_1x_4z_3r_3+z_1z_3z_4r_4+x_2y_2  \notag\\
&  & e^{szy}_{2}= -x_1x_3x_4r_1-x_3z_1z_4r_2-x_1z_3z_4r_2-x_4z_1z_3r_3+z_2y_2  \notag \\
&  & e^{sxy}_{3} = x_1x_2z_4r_2+x_2x_4z_1r_3+x_1x_4z_2r_3+z_1z_2z_4r_4+x_3y_3 \notag \\
&  & e^{szy}_{3} = -x_1x_2x_4r_1-x_2z_1z_4r_2-x_1z_2z_4r_2-x_4z_1z_2r_3+z_3y_3 \notag \\
&  & e^{sxy}_{4} = x_2x_3z_1r_2+x_1x_3z_2r_2+x_1x_2z_3r_2+z_1z_2z_3r_4+x_4y_4 \notag \\
&  & e^{szy}_{4} = -x_1x_2x_3r_1-x_3z_1z_2r_3-x_2z_1z_3r_3-x_1z_2z_3r_3+z_4y_4 \notag  \\
 &  &  e^{sy}_{12} \;\,  = -x_3^2z_4^2r_2^2+x_3^2x_4^2r_1r_3-x_3x_4z_3z_4r_2r_3-x_4^2z_3^2r_3^2+  \notag   \\
 &  &  \quad \quad \quad \quad   x_3x_4z_3z_4r_1r_4+z_3^2z_4^2r_2r_4+y_1y_2  \notag \\
 &  & e^{sy}_{13} \; \,= -x_2^2z_4^2r_2^2+x_2^2x_4^2r_1r_3-x_2x_4z_2z_4r_2r_3-x_4^2z_2^2r_3^2+ \label {eqn!explicits} \\
  &  &   \quad \quad \quad \quad      x_2x_4z_2z_4r_1r_4+ z_2^2z_4^2r_2r_4+y_1y_3  \notag  \\
&  & e^{sy}_{14}  \;\, = x_2^2x_3^2r_1r_2-x_3^2z_2^2r_2r_3-x_2x_3z_2z_3r_2r_3-x_2^2z_3^2r_2r_3+  \notag  \\
 &  &  \quad \quad \quad \quad   x_2x_3z_2z_3r_1r_4+z_2^2z_3^2r_3r_4+y_1y_4 \notag \\
&  & e^{sy}_{23} \;\, = -x_1^2z_4^2r_2^2+x_1^2x_4^2r_1r_3-x_1x_4z_1z_4r_2r_3-x_4^2z_1^2r_3^2+ \notag \\
  &  &  \quad \quad \quad \quad    x_1x_4z_1z_4r_1r_4+z_1^2z_4^2r_2r_4+y_2y_3 \notag \\
&  & e^{sy}_{24}  \; \, = x_1^2x_3^2r_1r_2-x_3^2z_1^2r_2r_3-x_1x_3z_1z_3r_2r_3-x_1^2z_3^2r_2r_3+ \notag   \\
  &  &   \quad \quad \quad \quad    x_1x_3z_1z_3r_1r_4+ z_1^2z_3^2r_3r_4+y_2y_4 \notag \\
 &  & e^{sy}_{34}  \;\, = x_1^2x_2^2r_1r_2-x_2^2z_1^2r_2r_3-x_1x_2z_1z_2r_2r_3-x_1^2z_2^2r_2r_3+  \notag  \\
   &  &  \quad \quad \quad \quad     x_1x_2z_1z_2r_1r_4+z_1^2z_2^2r_3r_4+y_3y_4  \notag   
\end {eqnarray}
By  (\ref{eqn!dfnofhi}) we have $ h_1 = z_1+z_2+z_3$, $h_2 = z_4$ and $h_3 = x_4 + r_5 (x_1+x_2+x_3)$.

We set $f^x_i \in k^{12}$, for $1 \leq i \leq 4$, to be the vector
with one on the coordinate corresponding to $x_i$ and
zero  elsewhere, $f^z_i \in k^{12}$  to be the vector
with one on the coordinate corresponding to $z_i$ and
zero  elsewhere,
and $ f^y_i \in k^{12}$ to be the vector
with one on the coordinate corresponding to $y_i$ and
zero  elsewhere.

Set $V_1$ to be the vector space spanned by $f^x_1,f^x_2,f^x_3$,
$V_2$ to be the vector space spanned by $f^z_1,f^z_2,f^z_3$,
$V_3$ to be the vector space spanned by $f^x_4$,
$V_4$ to be the vector space spanned by $f^z_4$,
$V_5$ to be the vector space spanned by $f^y_1,f^y_2,f^y_3$
and  $V_6$ to be the vector space spanned by $f^y_4$.
We have that, for $1 \leq i \leq 6$,  the vector space $V_i$ is $G$-invariant.
More precisely, using  (\ref{eqn!forgP}) we get 
\[
  gf^x_i =  -f^x_t, \quad \quad  gf^z_i =  f^z_t, \quad \quad  gf^y_i =  -f^y_t,
\]
where $t \in \{1,2,3 \} $ is uniquely  specified by $ t \equiv i+1 \mod 3$, and
\[
  gf^x_4 =  -f^x_4, \quad \quad  gf^z_4 =  f^z_4, \quad \quad  gf^y_4 =  -f^y_4.
\]

Since every element of $G$ different from the identity has a power
equal to $g^2$ or of $g^3$ to prove Proposition~3.18
it is enough to show
that if  $t \in \{ 2,3 \}$ and
\begin {equation}  \label{eqn!defnofP}
  P = \sum_{i=1}^4a^x_if^x_i + \sum_{i=1}^4a^z_if^z_i + \sum_{i=1}^4 a^y_if^y_i
 \in S^c_{cl},
\end {equation}  (with $a^x_i,a^z_i,a^y_i \in k$) are such that
there exists $h \in k^*$  with $g^t P = h*P$ then $P=0$, where by definition
\[
    h*P = \sum_{i=1}^4ha^x_if^x_i + \sum_{i=1}^4ha^z_if^z_i  +\sum_{i=1}^4h^2a^y_if^y_i
\]

\emph{Step 1.} We first study the action of $g^2$.  The action of $g^2$
on the direct sum of the $V_i$  can be described by
\[
    g^2  = (f^x_1,f^x_3,f^x_2) (f^z_1,f^z_3,f^z_2) (f^x_4) (f^z_4) (f^y_1,f^y_3,f^y_2)(f^y_4)
\]
in the sense that $g^2  f^x_1 = f^x_3, \;  \; g^2  f^x_3= f^x_2,  \;   \; g^2  f^x_2 = f^x_1 $ etc.
We assume
\begin {equation} \label {eqn!hypothesisforP}
    g^2(P) = h *P
\end{equation}
for some nonzero $P$  as in (\ref{eqn!defnofP}) and we will get a contradiction.
Since $P$ is nonzero, by  Lemma~3.8   there exists $i$,
with $1 \leq i \leq 4$, such that
$a^x_i \not= 0$ or $a^z_i \not= 0$.  Since the eigenvalues of $g^2$ acting on any
of $V_1, \dots ,V_4$ are contained in the set  $ \{ 1, \zeta^2, \zeta^4 \} $, we get that
$h \in \{ 1, \zeta^2, \zeta^4 \}$.

\emph{Step 2.} We assume that $h=1$  in (\ref{eqn!hypothesisforP}) and we will get a contradiction. By
looking at the action of $g^2$ on $V_1$ and $V_2$ we have
\[
   a^x_3 = a^x_2 = a^x_1, \quad   a^z_3 = a^z_2 = a^z_1. 
\]
Using  the equations $h_1,h_2,h_3$ we additionally get
\[
   a^z_1 = a^z_2 = a^z_3 = a^z_4 = 0, \quad   a^x_4 = -3r_5a^x_1.
\]
Substituting to  $ e^{szy}_{4}$ 
we get $(a^x_1)^3 r_1 =0$, hence $a^x_1 =0$ (since $r_1$ is
general), which implies that all $a^x_j =0$ and all $a^z_j =0$, contradicting Lemma~3.8.

\emph{Step 3.}  We assume that $h=\zeta^2$  in (\ref{eqn!hypothesisforP}) and we will get a contradiction. By
looking at the action of $g^2$ on each $V_i$  we have
\begin {eqnarray*}
  &  &  a^x_3 =  \zeta^2 a^x_1, \quad   a^x_2 = \zeta^4 a^x_1,  \quad   a^z_3 =\zeta^2 a^z_1,
         \quad  a^z_2 = \zeta^4 a^z_1,  \\
  &  & a^x_4  =  a^z_4 = a^y_4 = 0,  \quad  a^y_3 =\zeta^4 a^y_1, \quad a^y_2  =\zeta^2 a^y_1.
\end{eqnarray*}
Substituting to $ e^{sxy}_{1}$ 
we get $a^x_1a^y_1 =0$ and substituting to 
$  e^{szy}_{1}  $
we get $a^z_1a^y_1 =0$.
We can assume that $a^y_1 =0$, otherwise we get that all $a^x_j =0$ and all $a^z_j =0$ 
contradicting Lemma~3.8. As a consequence,  $a^y_i =0$ for $1 \leq i \leq 4$.

Substituting to $e^{sxy}_{4}$ we get 
\[
      3r_2(a^x_1)^2a^z_1 + r_4 (a^z_1)^3 =0,
\]
while substituting to $e^{szy}_{4}$ we get 
\[
    r_1(a^x_1)^3 + 3r_3a^x_1(a^z_1)^2 = 0.
\]
Hence 
\[
  r_1a^x_1e^{sxy}_{4} - 3r_2a^z_1e^{szy}_{4}  = (a^z_1)^3(a^x_1)(-9r_2r_3+r_1r_4),
\]
therefore $a^z_1 =0$ or $a^x_1 = 0$. If  $a^z_1 =0$ substituting to  $e^{szy}_{4}$ we get 
$a^x_1 = 0$, while if $a^x_1 =0$ substituting to  $e^{sxy}_{4}$ we get $a^z_1 =0$. 
In both cases all $a^x_j =0$ and all $a^z_j =0$  contradicting $P$ nonzero.

\emph{Step 4.}  We assume that $h=\zeta^4$  in (\ref{eqn!hypothesisforP}) and we will get a contradiction. By
looking at the action of $g^4$ on each $V_i$  we have
\begin {eqnarray*}
  &  &  a^x_3  =  \zeta^4 a^x_1, \quad   a^x_2 = \zeta^2 a^x_1,  \quad   a^z_3 =\zeta^4 a^z_1,
         \quad  a^z_2 = \zeta^2 a^z_1,  \\
  &  &  a^x_4  =  a^z_4 = a^y_4 = 0,  \quad  a^y_3 =\zeta^2 a^y_1, \quad a^y_2  =\zeta^4 a^y_1.
\end{eqnarray*}
Substituting to $e^{sxy}_{1}$ we get $a^x_1a^y_1 =0$ and substituting to $e^{szy}_{1}$ we get $a^z_1a^y_1 =0$.
We can assume that $a^y_1 =0$, otherwise we get that all $a^x_j =0$ and all $a^z_j =0$ 
contradicting Lemma~3.8. As a consequence,  $a^y_i =0$ for $1 \leq i \leq 4$.

Substituting to $e^{sxy}_{4}$ we get 
\[
      3r_2(a^x_1)^2a^z_1 + r_4 (a^z_1)^3 =0,
\]
while substituting to $e^{szy}_{4}$ we get 
\[
   3r_1(a^x_1)^3 + r_2a^x_1(a^z_1)^2 = 0.
\]
Hence 
\[
  0 =  r_1a^x_1e^{sxy}_{4} - 3r_2a^z_1e^{szy}_{4}  = (a^z_1)^3(a^x_1)(-9r_2r_3+r_1r_4),
\]
therefore $a^z_1 =0$ or $a^x_1 = 0$. If  $a^z_1 =0$ substituting to  $e^{szy}_{4}$ we get 
$a^x_1 = 0$, while if $a^x_1 =0$ substituting to  $e^{sxy}_{4}$ we get $a^z_1 =0$. 
In both cases all $a^x_j =0$ and all $a^z_j =0$  contradicting $P$ nonzero.

\emph{Step 5.} We now study the action of  $g^3$.
The action of $g^3$  on the  direct sum of the $V_i$  can be described by
\[
    g^3 (f^x_i)  =    - f^x_i, \quad \quad g^3 (f^z_i) =  f^z_i,  \quad \quad
    g^3 (f^y_i)  =    - f^y_i
\]
for $1 \leq i \leq 4$.

We assume
\begin {equation} \label {eqn!hypothesisforPforg3}
    g^3(P) = h *P
\end{equation}
for some nonzero $P$  as in (\ref{eqn!defnofP}) and we will get a contradiction.
Arguing as in Step 1, we get $h \in \{ 1, -1 \}$.

\emph{Step 6.}  We assume that $h=1$ in (\ref {eqn!hypothesisforPforg3})  and we will get a contradiction.
Indeed, by the way $g^3$ acts we have   $a^x_i = a^y_i =0$ for $1 \leq i \leq 4$.
Substituting to  $e^{sy}_{14}, e^{sy}_{24}, e^{sy}_{34}$ we get respectively
\[
      a^z_2a^z_3 = a^z_1a^z_3 = a^z_1a^z_2 = 0,
\]
hence at least two of  the three $a^z_1,a^z_2,a^z_3$ are zero. Using $h_1$ and $h_2$ we 
get that $a^z_i =0$ for $1 \leq i \leq 4$, contradicting $P$ nonzero.

\emph{Step 7.}  We assume that $h=-1$ in (\ref {eqn!hypothesisforPforg3})  and we will get a contradiction.
Indeed, by the way $g^3$ acts we have  $a^z_i =a^y_i = 0 $ for $1 \leq i \leq 4$.
Fix $1 \leq i < j \leq 4$. Substituting to $e^{sy}_{ij}$ we get that $a^x_a a^x_b = 0$
where $a,b$ have the property  $\{a,b,i,j \} = \{1,2,3,4 \}$.  As a consequence at least
three of the four $a^x_i$ are zero.  Using $h_3$ we get that all $a^x_i$ are zero contradicting 
$P$ nonzero. This finishes the proof of
Proposition~3.14.

\section {Final remarks and questions}   \label{sec!finalcom}

\begin{rem} \label {remark!no332}
In \cite{CR}, Corti and Reid pose the problem of interpreting
the Gorenstein formats arising from unprojection as
solutions to universal problems. What can be said about
the generic $\binom {n} {2}$ Pfaffians ideal of Definition~\ref{dfn!binomfaffiansideal}?
\end{rem}

\begin{rem} \label {remark!no333}
During the proof of Theorem~\ref{theorem!existenceCampedelli},
we established  that the \'{e}tale six to one numerical Campedelli covers
of our construction  have trivial algebraic fundamental group.
We expect that they also have trivial topological fundamental
group, but we were unable to prove it.
\end{rem}

\begin{rem} \label {remark!no334}
What is the dimension of the family of numerical Campedelli surfaces  $S/G$ of
Theorem~\ref{theorem!existenceCampedelli}?

We think it is unlikely, but we do not know it for certain, that the family of
surfaces $S/G$ of Theorem~\ref{theorem!existenceCampedelli}
gives the complete classification of numerical Campedelli surfaces with torsion $\Z/6$.
However, we expect that more refined geometric arguments
and unprojection machinery will eventually lead to a complete classification.
\end {rem}

\begin{rem} \label {remark!no335}
We believe that the ideas of the present paper can also be useful for the study of 
the numerical Campedelli surfaces with torsion groups $\Z/2$ and $\Z/3$.

Consider first the $\Z/3$ torsion case. The numeric invariants  suggest that the \'{e}tale
three to one  cover  of such a numerical Campedelli surface could be a suitable
member of $|-2K_{V_3}|$, where 
\[
    V_3 \subset \P(1^2,2^7,3^5)
\]
is a (candidate) codimension ten Fano threefold
having a basket of ten $1/2(1,1,1)$ singularities, which appears in
Brown's online database  of graded ring \cite{Br}. Moreover,
\cite{Br} suggests that $V_3$ can, perhaps, be constructed as a
result of a series of symmetric looking  type II unprojections
(cf.~\cite{R1}, \cite{P3}).

Similarly,  the numeric invariants for the $\Z/2$ torsion case
suggests that  the  \'{e}tale  two to one cover of such a numerical Campedelli could be
a suitable member of $|-2K_{V_2}|$, where 
\[
    V_2 \subset  \P(1,2^6,3^8)
\]
is a (candidate) codimension eleven Fano threefold
having a basket of twelve $1/2(1,1,1)$ singularities and also appearing
in \cite{Br}. Moreover, \cite{Br} suggests that $V_2$ can, perhaps,
be constructed as a  result of a series of again symmetric looking 
type IV unprojections  (cf.~\cite{R3}).

\end {rem}

\bigskip
\noindent
Jorge Neves \\
CMUC, Department of Mathematics, \\
University of Coimbra. \\
3001-454 Coimbra, Portugal. \\
e-mail: neves@mat.uc.pt

\bigskip
\noindent
Stavros Papadakis \\
Center for Mathematical Analysis, Geometry, and Dynamical Systems \\
Departamento de Matem\'atica,  Instituto Superior T\'ecnico \\
Universidade T\'ecnica de Lisboa \\
1049-001 Lisboa, Portugal \\
e-mail: spapad@maths.warwick.ac.uk

\end{document}